\documentclass[fleqn,twoside,twocolumn,nofootinbib]{revtex4} 
\usepackage{ujp} 
\begin{document}
\title[ON $T(n,4)$ TORUS KNOTS AND CHEBYSHEV POLYNOMIALS]
{ON \boldmath$T(n,4)$ TORUS KNOTS AND CHEBYSHEV POLYNOMIALS}%
\author{A.M. Pavlyuk}
\affiliation{Bogolyubov Institute for Theoretical Physics, Nat. Acad. of Sci. of Ukraine}
\address{14b, Metrolohichna Str., Kyiv 03680, Ukraine}
\email{pavlyuk@bitp.kiev.ua}

\udk{???} 

\razd{\secix}
\setcounter{page}{439}%
\maketitle

\newcommand{\ed}{\end{document}}
\newcommand{\be}{\begin{equation}}
\newcommand{\ee}{\end{equation}}
\newcommand{\bea}{\begin{eqnarray}}
\newcommand{\eea}{\end{eqnarray}}
\newcommand{\bc}{\begin{center}}
\newcommand{\ec}{\end{center}}
\newcommand{\ba}{\begin{array}}
\newcommand{\ea}{\end{array}}
\newcommand{\ve}{\varepsilon}
\newcommand{\wt}{\widetilde}
\newcommand{\cl}{\centerline}
\newcommand{\ts}{\textstyle}

\begin{abstract}
The Alexander polynomials $\Delta_{n,3}(t)$ and $\Delta_{n,4}(t)$ are presented as
a sum of the
 Alexander polynomials $\Delta_{k,2}(t).$
 These polynomials are also expressed in the form of a sum
 of Chebyshev polynomials of the second kind.
These expansions allow one to introduce the ``coordinates'' in
corresponding bases, which are proposed to be the numerical invariants
characterizing links and knots.
\end{abstract}

\section{Introduction}

The knot history began in 1867, when Lord Kelvin (W.H.~Thomson)
suggested to describe the atoms as knotted vortex tubes in the
ether~\cite{Th}. In 1975, L.D.~Faddeev proposed that knot-like
solitons could exist in a modified sigma model extended to the
three-dimensional space~\cite{Fa}. But the second advent of the
knots into physics~\cite{Wi,At,Ka,RV} began only in 1997 after the
article by Faddeev and Niemi~\cite{FN}. They made first attempts of
a numerical construction of the solitons with minimal energy in the
form of knots. Increasing the computer power demonstrated that a
number of linked and knotted configurations do exist in the Faddeev
model, which are the solutions characterized by local or global
energy minima.\looseness=1

Powerful tools in the knot theory are polynomial invariants.
In this paper, we concentrate on studying the Alexander polynomial
invariants for torus knots $T(n,l)\,,\,\,l{=}2,3,4\,,$
and on their connection with the Chebyshev polynomials~\cite{GP1} of the second kind.
The Alexander polynomials
$\Delta_{n,3}(t)$ and $\Delta_{n,4}(t)$ are presented firstly
as a sum of the
 Alexander polynomials $\Delta_{k,2}(t)$ and, secondly,
 through a sum
 of Chebyshev polynomials of the second kind.
These expansions allow us to introduce the ``coordinates'' in
corresponding bases, which are proposed, in the case of finding
similar expansions for all links and knots,
 to be numerical invariants characterizing
links and knots. The obtained expansions can be also used for the
investigation of baryon masses~\cite{Ga}.

\section{Alexander Polynomial Invariants for Torus Knots}

 The skein relation for the Alexander polynomials $\Delta(t)$ for knots and links,
  \be\label{a-skein}
\Delta_{+}(t)=(t^{1\over2}-t^{-{1\over2}})\Delta_{O}(t)+\Delta_{-}(t)\,
 \ee
 together with the condition for unknots,
  \be\label{unknot} \Delta_{\rm unknot}=1\,
  \ee
 gives an axiomatic definition of the Alexander polynomials.
 Using (\ref{a-skein}) and (\ref{unknot}), one can find the
Alexander polynomial for any knot or link 
with the help of the so-called surgery
 operations of ``elimination'' and ``switching''.
We denote torus knots as $T(n,l)$,  where  $n$ and $l$ are coprime
positive integers, and the corresponding  Alexander polynomial as $\Delta_{n,l}(t)$.

 As known~\cite{Ro, Li}, the Alexander polynomial  for the torus knot
$T(n,l)$ is given by the formula
 \be\label{a-nl}\ba{l}
 {\Delta}_{n,l}(t)=
{{(\ts t^{nl\over 2}-t^{-{nl\over
2}})(t^{1\over 2}-t^{-{1\over 2}})} \over{(\ts t^{n\over
2}{-}t^{-{n\over 2}})(t^{l\over 2}{-}t^{-{l\over 2}})}}\,
 \ea\ee
and has the form of the Laurent polynomial with the highest positive degree
\be\label{max}m={1\over2}(n-1)(l-1)\,.\ee

For $l=2\,,$ (\ref{a-nl}) gives
 \be\label{a-n2}
{\Delta}_{n,2}(t)= {{t^{n\over 2}+t^{-{n\over 2}}} \over{t^{1\over
2}+t^{-{1\over 2}}}}\,.
 \ee
Some first examples of $\Delta_{n,2}(t)$ from (\ref{a-n2})
are
 \bea \label{a-n2-ex}
\hspace{-0.5cm}&&{\Delta}_{1,2}(t)=1,\quad
{\Delta}_{3,2}(t)=t-1+t^{-1}, \nonumber\\\vspace{2mm}
\hspace{-0.5cm}&&{\Delta}_{5,2}(t)=t^{2}-t+1-t^{-1}+t^{-2}
\,. \eea  Below, we will find an expansion of the Alexander
polynomials for two sets of torus knots, $T(n,3)$ and $T(n,4)\,,$
through the Alexander polynomials for torus knots $T(k,2)$.

If $l=3\,,$ (\ref{a-nl}) looks as
 \be\label{a-n3} {\Delta}_{n,3}(t)= {{t^{n}+1+t^{-{n}}}
\over{t^{1}+1+t^{-{1}}}}\,.
 \ee
Some  examples of (\ref{a-n3}) are
 \bea
 \label{a-n3-ex}
\hspace{-0.5cm}&&{\Delta}_{1,3}(t)=1,\quad
{\Delta}_{2,3}(t)=t-1+t^{-1}, \nonumber\\\vspace{2mm}
\hspace{-0.5cm}&&{\Delta}_{4,3}(t)=t^{3}-t^{2}+1-t^{-2}+t^{-3},
\nonumber\\\vspace{2mm}
\hspace{-0.5cm}&&{\Delta}_{5,3}(t)=t^{4}-t^{3}+t-1+t^{-1}-t^{-3}+t^{-4}\,.
 \eea

For $l=4\,,$ relation (\ref{a-nl}) yields
 \be\label{a-n4}
 {\Delta}_{n,4}(t)= {t^{\frac{3n}{2}}+t^{\frac{n}{2}}+t^{-{\frac{n}{2}}}+t^{-{\frac{3n}{2}}}
\over{t^{\frac{3}{2}}+t^{\frac{1}{2}}+t^{-{\frac{1}{2}}}+t^{-{\frac{3}{2}}}}}\,.
 \ee
The first examples of (\ref{a-n4}) are as follows:
 \bea
 \label{a-n4-ex}
\hspace{-0.5cm}&&{\Delta}_{1,4}(t)=1,\quad
{\Delta}_{3,4}(t)=t^{3}-t^{2}+1-t^{-2}+t^{-3},
\nonumber\\\vspace{2mm}
\hspace{-0.5cm}&&{\Delta}_{5,4}(t)=t^{6}-t^{5}+t^{2}-1+t^{-2}+t^{-5}+t^{-6}
\,.
 \eea

\section{On the Expression of \boldmath$\Delta_{n, 3}(t)$ through \boldmath$\Delta_{k, 2}(t)$}

In this section, we find the formula expressing any Alexander polynomial
$\Delta_{n, 3}(t)$ as an algebraic sum of
the Alexander polynomials $\Delta_{k,2}(t)$~\cite{GP2}.

{\bf{Proposition 1.}}
{\it The Alexander polynomial for any torus knot $\, T(n,3)$
is presented as an algebraic sum of the Alexander polynomials
for torus knots $\, T(k,2)$ with the help of the formula
 \be\label{p1} {\Delta}_{n,3}(t)
=\sum\limits_{j=0}^{d_{1}}
{\Delta}_{2n-1-6j,\, 2}(t)-
\sum\limits_{i=0}^{d_{2}}
{\Delta}_{2n-5-6i,\, 2}(t)\bigr)\,,
  \ee
where
$d_{1}=\Bigl[\frac{2n-1}{6}\Bigr]$ and
$d_{2}=\Bigl[\frac{2n-5}{6}\Bigr]$
}.

The proof of (\ref{p1}) follows from the relation
\be\label{p1-}
{\Delta}_{n,3}(t)-{\Delta}_{n-3,3}(t)={\Delta}_{2n-1,2}(t)-{\Delta}_{2n-5,2}(t)\,,\ee
which can be proved, in turn, with the help of (\ref{a-n2}) and (\ref{a-n3}).

Equation (\ref{p1}) can be rewritten as
  \be\ba{l}\label{t2-}
 {\Delta}_{n,3}(t)\equiv {\Delta}_{3d+r,3}(t)=
\vspace{2mm}\\
=\sum\limits_{j=0}^{d}
{\Delta}_{2n-1-6j,2}(t)-
\sum\limits_{i=0}^{d-1}
{\Delta}_{2n-5-6i,2}(t)
\,,
 \ea\ee
where $r=1,2,$ and $d=0,1,2,3\ldots\ .$ For $n$ given, $d$ is
the integer part of $\frac{n}{3}:$
$\,\,\,d=\Bigl[\frac{n}{3}\Bigr].$

Using (\ref{p1}), we give a few examples expressing the
Alexander polynomials ${\Delta}_{n,3}(t)$ as a sum of the
Alexander polynomials ${\Delta}_{k,2}(t):$
 \bea  \label{exam32}
\hspace{-0.5cm}&&{\Delta}_{1,3}(t)={\Delta}_{1,2}(t)\,,\quad
{\Delta}_{2,3}(t)={\Delta}_{3,2}(t)\,,\quad \nonumber\\\vspace{2mm}
\hspace{-0.5cm}&&{\Delta}_{4,3}(t)={\Delta}_{7,2}(t)-{\Delta}_{3,2}(t)+{\Delta}_{1,2}(t)\,,\quad
\nonumber\\\vspace{2mm}
\hspace{-0.5cm}&&{\Delta}_{5,3}(t)={\Delta}_{9,2}(t)-{\Delta}_{5,2}(t)+{\Delta}_{3,2}(t)\,,\quad
\nonumber\\\vspace{2mm}
\hspace{-0.5cm}&&{\Delta}_{7,3}(t)={\Delta}_{13,2}(t)-{\Delta}_{9,2}(t)+{\Delta}_{7,2}(t)-
\nonumber\\\vspace{2mm}
\hspace{-0.5cm}&&-{\Delta}_{3,2}(t)+{\Delta}_{1,2}(t)\,,
\nonumber\\\vspace{2mm}
\hspace{-0.5cm}&&{\Delta}_{8,3}(t)={\Delta}_{15,2}(t)-{\Delta}_{11,2}(t)+{\Delta}_{9,2}(t)-
\nonumber\\\vspace{2mm}
\hspace{-0.5cm}&&-{\Delta}_{5,2}(t)+{\Delta}_{3,2}(t)
\,.
 \eea
 It is easy to see that the expansion of ${\Delta}_{n,3}(t)$
 includes $2d+1$ terms like ${\Delta}_{k,2}(t)$.

\section{On ``Coordinates''\ of \boldmath$\Delta_{n, 3}(t)$ in the \boldmath$\Delta_{k, 2}(t)$-Basis }

From (\ref{gavr}), it follows that the Alexander polynomials are
orthonormal due to the orthonormality of the Chebyshev polynomials (with
a definite weight on the interval $[-1,1]$).

Thus, we can introduce an orthonormal basis consisting of
$\Delta_{k,2}(t)\,$ and find an expansion of any $\Delta_{n,3}(t)$
in this basis. We may consider $k$ as the
``coordinates'' of a knot $T(n,3)$ in the $\Delta_{k,2}(t)$-basis.
From (\ref{exam32}), we have the presentation of the particular torus
knots using such ``coordinates'':
\bea \label{coord3}  \hspace{-0.5cm}&&{\Delta}_{1,3}(t){=}(1),\quad
{\Delta}_{2,3}(t){=}(3),\quad {\Delta}_{4,3}(t){=}(7; -3; 1),
\nonumber\\\vspace{2mm} \hspace{-0.5cm}&&{\Delta}_{5,3}(t){=}(9; -5;
3),\quad {\Delta}_{7,3}(t){=}(13; -9; 7; -3; 1),
\nonumber\\\vspace{2mm} \hspace{-0.5cm}&&{\Delta}_{8,3}(t){=}(15;
-11; 9; -5; 3), \nonumber\\\vspace{2mm}
\hspace{-0.5cm}&&{\Delta}_{10,3}(t){=}(19; -15; 13; -9; 7; -3; 1)\,.
\eea
 Could we find a similar representation for an arbitrary knot, we gain the
numerical invariant for knots (and links).

\section{On the Expression of \boldmath$\Delta_{n, 4}(t)$ through \boldmath$\Delta_{k, 2}(t)$}

In this section, we present any Alexander polynomial
$\Delta_{n, 4}(t)$ as an algebraic sum of Alexander polynomials $\Delta_{k,2}(t)$
of $T(k,2)$ torus knots.

{\bf{Proposition 2.}} 
{\it The Alexander polynomial for any torus knot $\, T(n,4)$
is presented as a sum of the Alexander polynomials
for torus knots $\, T(k,2)$ by the formula
 \be\ba{l}\label{p2}
 {\Delta}_{n,4}(t)=\sum\limits_{j=0}^{d_1}(-1)^{j}
 {\Delta}_{3n-2-4j,2}(t)+
\vspace{2mm}\\
+ \sum\limits_{i=0}^{d_2}(-1)^{i} {\Delta}_{n-2-4i,2}(t) \,,
 \ea\ee
where
$\,d_{1}=\Bigl[\frac{3n-2}{4}\Bigr]$
and
$\,d_{2}=\Bigl[\frac{n-2}{4}\Bigr].$
}

To prove (\ref{p2}), we use
\bea\label{p2-}
\hspace{-0.5cm}&&{\Delta}_{n,4}(t)+{\Delta}_{n-4,4}(t)=
{\Delta}_{3n-2,2}(t)-{\Delta}_{3n-6,2}(t)+
\nonumber\\\vspace{2mm}
\hspace{-0.5cm}&&+{\Delta}_{3n-10,2}(t)+{\Delta}_{n-2,2}(t) , \eea
which can be proved, in turn, with the help of (\ref{a-n3}) and
(\ref{a-n4}).

Using (\ref{p2}), we give some first examples of the presentation of the
Alexander polynomials ${\Delta}_{n,4}(t)$ by the a of the
Alexander polynomials ${\Delta}_{k,2}(t):$
 \bea \label{exam42}
\hspace{-0.5cm}&&{\Delta}_{1,4}(t)={\Delta}_{1,2}(t),
\nonumber\\\vspace{2mm}
 \hspace{-0.5cm}&&{\Delta}_{3,4}(t)=
{\Delta}_{7,2}(t)-{\Delta}_{3,2}(t)+{\Delta}_{1,2}(t),
\nonumber\\\vspace{2mm}
 \hspace{-0.5cm}&&{\Delta}_{5,4}(t)=
{\Delta}_{13,2}(t)-{\Delta}_{9,2}(t)+{\Delta}_{5,2}(t)+
\nonumber\\\vspace{2mm}
 \hspace{-0.5cm}&&+{\Delta}_{3,2}(t)-{\Delta}_{1,2}(t),
 \nonumber\\\vspace{2mm}
\hspace{-0.5cm}&&{\Delta}_{7,4}(t)={\Delta}_{19,2}(t)-{\Delta}_{15,2}(t)+{\Delta}_{11,2}(t)-
\nonumber\\\vspace{2mm}
 \hspace{-0.5cm}&&-{\Delta}_{7,2}(t)
+{\Delta}_{5,2}(t) +{\Delta}_{3,2}(t)-{\Delta}_{1,2}(t) .
\eea
 The expansion of ${\Delta}_{n,4}(t)$ includes $n$ terms like ${\Delta}_{k,2}(t)$.

Thus, relation (\ref{exam42}) yields the ``coordinates'' of
corresponding $T(n,4)$ knots (i.e., the ``coordinates'' of a knot
$T(n,4)$ in the $\Delta_{k,2}(t)$-basis)
 \bea \label{coord3}
    \hspace{-0.5cm}&&{\Delta}_{1,4}(t){=}(1),\quad
{\Delta}_{3,4}(t){=}(7; -3; 1), \nonumber\\\vspace{2mm}
 \hspace{-0.5cm}&&{\Delta}_{5,4}(t){=}(13; -9; 5; 3; -1), \nonumber\\\vspace{2mm}
 \hspace{-0.5cm}&&{\Delta}_{7,4}(t){=}(19; -15; 11; -7; 5; 3; -1),
\nonumber\\\vspace{2mm}  \hspace{-0.5cm}&&{\Delta}_{9,4}(t){=}(25;
-21; 17; -13; 9; 7; -5; -3; 1) . \eea

In the following section, we are going to present the Alexander polynomials
$\Delta_{n,2}(t)$,\ $\Delta_{n,3}(t),$ and $\Delta_{n,4}(t)$ in terms of the Chebyshev
polynomials of the second kind $V_{n}(x)\,,\,\, x=t+t^{-1}\,. $

\section{Alexander Polynomials in Terms of Chebyshev Polynomials}

The monic Chebyshev polynomials of the second kind,
 which have unit coefficients of
$x^{n}$, can be defined as
   \be\label{V}
   V_{n}(x)={\sin(n+1)\theta\over {\sin\theta}}\,,\quad
2\cos\theta=x
  \ee
 or by the recurrence relation
 \be\label{rec-V} V_{n+1}=xV_{n}-V_{n-1}\,,\quad V_{0}=1\,,\quad
V_{1}=x\,.
 \ee
 Some first low-order Chebyshev polynomials of
the second kind are:
 \bea \label{Vx}
\hspace{-0.5cm}&& V_{0}=1\,,\quad V_{1}=x\,,\quad
V_{2}=x^{2}-1\,,\quad V_{3}{=}x^{3}{-}2x,
 \nonumber \\\vspace{2mm}
\hspace{-0.5cm}&&V_{4}=x^{4}-3x^{2}+1\,,\quad
V_{5}{=}x^{5}{-}4x^{3}{+}3x.
 \eea

By analogy to the ``coordinates'' of a knot 
in the $\Delta_{k,2}(t)$-basis, we introduce the ``coordinates'' of
a knot in the $\wt{V}_{k}(x)$-basis, where $\wt{V}_{k}(x)$ is a
certain combination of Chebyshev polynomials of the second kind,
namely
\[\wt{V}_{k}(x)=V_{k}(x)-V_{k-1}(x)-V_{k-2}(x)+
\]
\[{+}V_{k-3}(x){+}V_{k-4}(x){-}V_{k-5}(x){-}\cdots V_{0}(x){=} \]
\bea
\label{wtV}=\sum\limits_{j=1}^{k+1}(-1)^{[\frac{j}{2}]}V_{k+1-j}(x).
 \eea
For example, the $\wt{V}_{7}(x)$-basis
looks as
\[ \wt{V}_{7}(x)=V_{7}(x)-V_{6}(x)-V_{5}(x)+V_{4}(x)+V_{3}(x)- \]
\[ -V_{2}(x)-V_{1}(x)+V_{0}(x).\]

The Alexander polynomials $\Delta_{n,2}(t)$ satisfy relation~\cite{Ga,GP2}
\[
\Delta_{n,2}(t)\equiv\Delta_{2m+1,2}(t)=V_{m}(x)-V_{m-1}(x)\,,
\]
\bea\label{gavr}x=t+t^{-1}\,.
 \eea
From (\ref{wtV}) and (\ref{gavr}), we have
\be \Delta_{n,2}(t)\equiv\Delta_{2m+1,2}(t){=}{\wt V}_{m}(x){+}{\wt
V}_{m-2}(x), \ x{=}t{+}t^{-1}.
\ee
Whence, we obtain the
``coordinates'' of $\Delta_{n,2}(t)$\ in\ ${\wt V}$-basis:
\be\label{coordV2}
 \Delta_{2m+1,2}(t)=(m;\ m{-}2)_{\wt V}\,,
\ee
which can be written as
\be\label{coordV2-} \Delta_{n,2}(t)= {\wt V}_{\frac{n-1}{2}}+{\wt
V}_{\frac{n-5}{2}}= \Bigl(\frac{n-1}{2};\ \frac{n-5}{2}\Bigr)_{\wt
V}\,.\ee

Let us write the  ``V-coordinates'' for  $\Delta_{n,3}(t)$. It follows from
(\ref{exam32}) and (\ref{gavr}) that
\[
 {\Delta}_{n,3}(t)
=V_{n-1}(x)+
\]
\[{+}\sum\limits_{k=0}^{d}\bigl({-}V_{n-2-3k}(x){-}V_{n-3-3k}(x){+}
2V_{n-4-3k}(x)\bigr){=}
\]
\be\label{n3V}=\sum\limits_{j=0}^{d_{1}}{\wt
V}_{n-1-3j}(x)-\sum\limits_{i=0}^{d_{2}}{\wt V}_{n-5-3i}(x) \,, \ee
where $\ d{=}\Bigl[\frac{n-2}{3}\Bigr], \
d_{1}{=}\Bigl[\frac{n-1}{3}\Bigr], \
d_{2}{=}\Bigl[\frac{n-5}{3}\Bigr], \
x{=}t+t^{-1}.$\\
Below, we give some examples of (\ref{n3V}):
 \bea\label{n3V-}
\hspace{-0.5cm}&&\Delta_{1,3}(t)=V_{0}(x)={\wt V}_{0}(x)=(0)_{\wt
V}, \nonumber\\\vspace{2mm} \hspace{-0.5cm}
\hspace{-0.5cm}&&\Delta_{2,3}(t)=V_{1}(x)-V_{0}(x)={\wt
V}_{1}(x)=(1)_{\wt V}, \nonumber\\\vspace{2mm}
\hspace{-0.5cm}&&\Delta_{4,3}(t)=V_{3}(x)-V_{2}(x)-V_{1}(x)+2V_{0}(x)=
\nonumber\\\vspace{2mm}
\hspace{-0.5cm}&& ={\wt V}_{3}(x)-{\wt
V}_{0}(x)=(3;\, -0)_{\wt V}, \nonumber\\\vspace{2mm}
\hspace{-0.5cm}&&\Delta_{5,3}(t)=V_{4}(x){-}V_{3}(x){-}V_{2}(x){+}2V_{1}(x){-}V_{0}(x){=}
\nonumber\\\vspace{2mm}
\hspace{-0.5cm}&& ={\wt V}_{4}(x)+{\wt
V}_{1}(x)-{\wt V}_{0}(x)=(4;\, 1; \, -0)_{\wt V} \,.
    \eea

{\bf{Proposition 3.}}
{\it The Alexander polynomial for the torus knot $T(n,4)$ is expressed
as a sum of Chebyshev polynomials of the second kind by the
relation
 \bea\label{p3} \hspace{-0.5cm}&&{\Delta}_{n,4}(t)
= \nonumber\\\vspace{2mm}
\hspace{-0.5cm}&&{=}\sum\limits_{j=1}^{\frac{3n-1}{2}}(-1)^{[\frac{j}{2}]}V_{\frac{3n-1}{2}-j}(x){+}
\sum\limits_{i=1}^{\frac{n-1}{2}}({-}1)^{[\frac{i}{2}]}V_{\frac{n-1}{2}-i}(x){=}
\nonumber\\\vspace{2mm} \hspace{-0.5cm}&&={\wt
V}_{\frac{3n-3}{2}}+{\wt V}_{\frac{n-3}{2}}= \Bigl(\frac{3n-3}{2};\
\frac{n-3}{2}\Bigr)_{\wt V}, \eea
where $ x=t+t^{-1}\,.$}

Some first examples of (\ref{p3}) are as follows:
\bea \label{p3a} \hspace{-0.5cm}&&\Delta_{1,4}(t)=V_{0}(x)={\wt
V}_{0}(x)=(0)_{\wt V}, \nonumber\\\vspace{2mm}
\hspace{-0.5cm}&&\Delta_{3,4}(t)=V_{3}(x)-V_{2}(x)-V_{1}(x)+2V_{0}(x)=
\nonumber\\\vspace{2mm} \hspace{-0.5cm}&&
 ={\wt V}_{3}(x)+{\wt
V}_{0}(x)=(3;\,0)_{\wt V}, \nonumber\\\vspace{2mm}
\hspace{-0.5cm}&&\Delta_{5,4}(t)=V_{6}(x)-V_{5}(x)-V_{4}(x)+V_{3}(x)+
\nonumber\\\vspace{2mm}
 \hspace{-0.5cm}&&+V_{2}(x){-}2V_{0}(x){=}{\wt
V}_{6}(x){+}{\wt V}_{1}(x){=}(6;\,1)_{\wt V} \,.
 \eea

\section{Concluding Remarks}

Here, we have shown how to introduce the ``coordinates'' for
torus knots and links with the help of corresponding Alexander
polynomials. As a basis, we can choose the simplest
Alexander polynomials $\Delta_{k,2}(t)$ or the Chebyshev polynomials
of the second kind. We have found the above-mentioned expansions for
the torus knots $\Delta_{n,l}(t),\ l{=}2,3,4.$ By finding
such expansions for all (torus) knots and links, we have obtained the numerical
invariants for them.

\rezume{%
ПРО ТОРИЧНІ ВУЗЛИ $T(n,4)$ І ПОЛІНОМИ ЧЕБИШОВА}
{%
А.М. Павлюк} {Поліноми Александера $\Delta_{n,3}(t)$ і
$\Delta_{n,4}(t)$ представлено як суму поліномів Александера
$\Delta_{k,2}(t)$. Ці поліноми також виражено через суму поліномів
Чебишова другого роду. Отримані розклади дозволяють ввести
``координати'' щодо вказаних базисів, які, як передбачається, є
числовими інваріантами вузлів і зачеплень.}


\begin{thebibliography}{99}
\bibitem{Th}  W.H. Thomson, Proc. of Roy. Soc. Edinburg
 {\bf 6}, 94 (1867).

 \bibitem{Fa}  L.D. Faddeev, {\it Quantization of Solitons},
Preprint IAS-75-QS70 (Inst. for Advanced Study, Princeton, 1975).

\bibitem{Wi} E. Witten, Comm. Math. Phys. {\bf 121}, 351 (1989).

\bibitem{At} M.F. Atiyah, {\it The Geometry and Physics of Knots}  (Cambridge Univ. Press, Cambridge, 1990).

\bibitem{Ka} L.H. Kauffman, {\it Knots and Physics}  (World Sci., Singapore, 2001).

\bibitem{RV} E. Radu and M.S. Volkov,
 Phys. Rep. {\bf 468}, No.~4, 101 (2008); arXiv:0804.1357v2 [hep-th].

\bibitem{FN}  L. Faddeev and A.J. Niemi,
 Nature {\bf 387}, 58 (1997); arXiv:hep-th/9610193.

\bibitem{GP1} A.M. Gavrilik and A.M. Pavlyuk,
  Ukr. J. Phys. {\bf 55}, 129 (2010);
  arXiv:0912.4674v2 [math-ph].

\bibitem{Ga} A.M. Gavrilik,
 J. Phys. A {\bf 27}, 91 (1994);
Nucl. Phys. B (Proc. Suppl.)  {\bf 102}, 298 (2001),
arXiv:hep-th/0103325v4.

\bibitem{Ro} D. Rolfsen, {\it Knots and Links} (Amer. Math. Soc., Providence, RI, 2003).

\bibitem{Li} W.B.R. Lickorish, {\it An Introduction to Knot Theory} (Springer, New York, 1997).

\bibitem{GP2} A.M. Gavrilik and A.M. Pavlyuk,
  Ukr. J. Phys. {\bf 56}, 680 (2011);  arXiv:1107.5516v1 [math-ph].




\begin{flushright}
{\footnotesize Received 15.07.2011}
\end{flushright}
\end{thebibliography}
\end{document}